\documentclass[11pt, a4paper, reqno]{amsart}
\usepackage{preambolo}

\begin{document}


\setcounter{secnumdepth}{3}

\setcounter{tocdepth}{2}

\newcommand\hidemath{$(d + d^c)$}

\title{\textbf{On the spaces of $(d + d^c)$-harmonic forms and $(d + d^\Lambda)$-harmonic forms on almost Hermitian manifolds and complex surfaces}
}

\author[Lorenzo Sillari]{Lorenzo Sillari}

\address{Lorenzo Sillari: Scuola Internazionale Superiore di Studi Avanzati (SISSA), Via Bonomea 265, 34136 Trieste, Italy.} 
\email{lsillari@sissa.it}

\author[Adriano Tomassini]{Adriano Tomassini}

\address{Adriano Tomassini: Dipartimento di Scienze Matematiche, Fisiche e Informatiche, Unità di Matematica e
Informatica, Università degli Studi di Parma, Parco Area delle Scienze 53/A, 43124, Parma, Italy}
\email{adriano.tomassini@unipr.it}

\maketitle

\begin{abstract} 
\noindent \textsc{Abstract}. We study the spaces of $(d+d^c)$-harmonic forms and $(d+d^\Lambda)$-harmonic forms, a natural generalization of the spaces of Bott-Chern harmonic forms, resp.\ sympelctic harmonic forms, from complex, resp.\ symplectic, manifolds to almost Hermitian manifolds. We apply the same techniques to compact complex surfaces, computing their Bott-Chern and Aeppli numbers and their spaces of $(d+d^\Lambda)$-harmonic forms. We give several applications to compact quotients of Lie groups by a lattice.

\end{abstract}

\blfootnote{  \hspace{-0.55cm} 
{\scriptsize 2020 \textit{Mathematics Subject Classification}. Primary: 32Q60, 53D15; Secondary: 53C15, 32J15.

\noindent The authors are partially supported by GNSAGA of INdAM. The second author is partially supported by the Project PRIN 2017 “Real and Complex Manifolds: Topology, Geometry and
holomorphic dynamics” (code 2017JZ2SW5).}}

\section{Introduction}\label{sec:intro}

Let $(M,J)$ be a compact complex manifold and set $d^c := J^{-1} d J$. Fix a $J$-compatible Hermitian metric and consider the \emph{Dolbeault Laplacian}
\[
\Delta_{\bar \partial} := \bar \partial \bar \partial^* + \bar \partial^* \bar \partial,
\]
the \emph{Bott-Chern Laplacian}
\[
\Delta_{d + d^c} := d d^c ( d d^c)^* + ( d d^c)^* d d^c + d^* d^c ( d^* d^c )^* + ( d^* d^c )^* d^* d^c + d^* d + (d^c)^* d^c
\]
and the \emph{Aeppli Laplacian}
\[
\Delta_{d d^c} := d d^c ( d d^c)^* + ( d d^c)^* d d^c +  d (d^c)^* ( d (d^c)^* )^* + ( d (d^c)^* )^* d (d^c)^* + d d^* + d^c (d^c)^*
\]
or, equivalently, the Laplacians $\Delta_{\partial + \bar \partial}$, $\Delta_{\partial \bar \partial}$, obtained replacing $d$,$d^c$ by $\partial$,$\bar \partial$ (the use of our notation instead of the usual $\Delta_{\bar \partial}$, $\Delta_{BC}$, $\Delta_A$ is to have a notation uniform with the one adopted for almost Hermitian manifolds). Denote by $\H^{p,q}_{\bar \partial}$, $\H^{p,q}_{d + d^c}$, $\H^{p,q}_{d d^c}$ their respective kernels computed on $(p,q)$-forms. By Hodge theory \cite{Sch07} and the compactness assumption, they are finite-dimensional vector spaces over $\C$, whose respective dimensions $h^{p,q}_{\bar \partial}$, $h^{p,q}_{d + d^c}$, $h^{p,q}_{d d^c}$, which in principle depend on the choice of metric, are actually invariants of the complex structure thanks to isomorphisms with Dolbeault, Bott-Chern and Aeppli cohomologies
\[
H^{p,q}_{\bar \partial} \cong \H^{p,q}_{\bar \partial}, \quad H^{p,q}_{BC} \cong \H^{p,q}_{d + d^c}, \quad H^{p,q}_{A} \cong \H^{p,q}_{d d^c}.
\]
A similar theory has been developed by Tseng and Yau \cite{TY12a, TY12b} for symplectic manifolds. If $(M,\omega)$ is a compact symplectic manifold, let $d^\Lambda$ be the symplectic adjoint of $d$. Then, fixed a compatible metric, consider the sympelctic Laplacians $\Delta_{d+d^\Lambda}$, $\Delta_{d d^\Lambda}$. Their kernels $\H^k_{d+d^\Lambda}$, $\H^k_{d d^\Lambda}$, computed on complex $k$-forms, are finite-dimensional vector spaces over $\C$ of dimensions $h^k_{d+d^\Lambda}$, $h^k_{d d^\Lambda}$ (cf.\ remark \ref{rem:ty}). The importance of the numbers $h^k_{d+d^c}$, resp.\ $h^k_{d+d^\Lambda}$, in the complex, resp.\ symplectic, case is paramount and well-established. For instance, they satisfy Fr\"olicher-type inequalities and characterize the $d d^c$-lemma, resp.\ $ d d^\Lambda$-lemma \cite{AT15a, AT13}.
\vspace{.2cm}

Let now $(M,J,\omega,g)$ be a compact almost Hermitian manifold, that is a manifold endowed with an \emph{almost} complex structure $J$, a $J$-compatible metric $g$, and a non-degenerate $2$-form 
$\omega$ (not necessarily $d$-closed) such that $\omega(J \cdot, \, \cdot) = g (\cdot, \, \cdot)$. We refer to $\omega$ as an \emph{almost symplectic structure}.\\
Let $P \in \{ \bar \partial, \partial + \bar \partial, \partial \bar \partial, d + d^c, d d^c, d + d^\Lambda, d d^\Lambda \}$. As observed by Hirzebruch \cite{Hir54} for $P= \bar \partial$, by Piovani and the second author \cite{PT22} for $P= \partial + \bar \partial, \partial \bar \partial$, by the authors of this paper in \cite{ST23b} for $P =  d + d^c, d d^c$, and in proposition \ref{prop:ellipticity} for $P = d + d^\Lambda, d d^\Lambda$, the Laplacian $\Delta_P$ is still elliptic and self-adjoint, even without assuming that $J$ is a complex structure or $d \omega =0$. Hence, its kernel computed on $(p,q)$ or $k$-forms, $\H^{p,q}_P$, $\H^k_P$, the so-called \emph{space of $P$-harmonic forms}, is a finite-dimensional vector space over $\C$ of dimension $h^{p,q}_P$, $h^k_P$. In view of the isomorphism valid in the complex case, Kodaira and Spencer posed the following problem for $h^{p,q}_{\bar \partial}$. 
\vspace{.2cm}

\textbf{Kodaira-Spencer's Problem.} (Kodaira-Spencer, see \cite[Problem 20]{Hir54})\\
Is $h^{p,q}_{\bar \partial}$ independent of the choice of (almost) Hermitian structure? If not, give some other definition of the $h^{p,q}_{\bar \partial}$ which depends only on the almost complex structure and which generalizes the $h^{p,q}_{\bar \partial}$ of a complex manifold.
\vspace{.2cm}

A negative answer was given by Holt and Zhang almost 70 years later \cite{HZ22}. Since $h^{p,q}_{\bar \partial}$ depends on the choice of metric, we ask the following question, which appears as a natural extension of Kodaira-Spencer's problem.

\begin{displayquote}
Is $h^{p,q}_P$/$h^k_P$ independent of the choice of almost Hermitian metric?
\end{displayquote}

A consistent part of the present paper is devoted to investigating the metric-independence of $h^k_P$ when $P = d+d^c, dd^c, d+d^\Lambda, d d^\Lambda$.
\vspace{.2cm}

This paper consists of two main parts. In the first part (sections \ref{sec:laplacians}, \ref{sec:ddc}, \ref{sec:ddlambda}, and \ref{sec:relations}), we study the spaces of $(d+d^c)$-harmonic and $(d+d^\Lambda)$-harmonic forms on compact almost Hermitian manifolds. We are mainly interested in understanding the properties of $\H^k_{d+d^c}$ and $\H^k_{d+d^\Lambda}$, in determining to which extent the numbers $h^k_{d+d^c}$, $h^k_{d+d^\Lambda}$ are independent of the choice of metric and thus define an almost complex, resp.\ almost symplectic, invariant, and in establishing relations among $h^k_{d+d^c}$, $h^k_{d+d^\Lambda}$ and topological invariants, like the Betti numbers.\\
For the basic properties of $\H^k_{d+d^c}$, we refer the reader to \cite{ST23b} (see also section \ref{sec:ddc}). Here we show that similar properties hold also for $\H^k_{d+d^\Lambda}$. For what concerns the metric independence of $h^k_{d+d^c}$, $h^k_{d+d^\Lambda}$, the most significant results are obtained in dimension $4$ for $k=2$, where we prove two direct sum decompositions (Theorems \ref{thm:ddc:dec} and \ref{thm:ddlambda:dec}).

\begin{thmx}
    Let $(M,J,\omega,g)$ be an almost Hermitian $4$-manifold and let $(\tilde{\omega},\tilde{g})$ be a Gauduchon metric in the same conformal class of $(\omega,g)$. Then
    \begin{itemize}
        \item for a suitable $2$-form $\gamma_0$, there is a decomposition
        \[
        \H^2_{d + d^c} = \C \langle \Tilde{\omega} + \gamma_0 \rangle \oplus \H^-_g \oplus \H^{(2,0)(0,2)}_J.
        \]
        In particular, $h^2_{d+d^c} = b^- +1 + h^-_J$ and it is metric independent;

        \item there is a decomposition
        \[
        \H^2_{d + d^\Lambda} = 
        \begin{cases}
        \C \langle \omega \rangle \oplus 
        \mathcal{PH}^2_{d + d^\Lambda} & \text{if } d \omega = 0,\\
        \mathcal{PH}^2_{d + d^\Lambda} & \text{if } d \omega \neq 0,
        \end{cases}
        \]
        where $\mathcal{PH}^2_{d + d^\Lambda}$ denotes the space of $(d + d^\Lambda)$-harmonic primitive $2$-forms.
    \end{itemize}
\end{thmx}

Assuming that either $J$ is complex or $d \omega =0$, we can establish strong relations between $h^k_{d+d^c}$, $h^k_{d+d^\Lambda}$ (Theorems \ref{thm:ak:betti}, \ref{thm:hermitian:betti}, and \ref{thm:akinclusion}).

\begin{thmx}

Let $(M,J, \omega,g)$ be a Hermitian $2m$-manifold. Suppose that we have $\H^k_{d^\Lambda} \subseteq \H^k_{d + d^\Lambda}$ for some $k$. Then $\H^k_{d + d^\Lambda} = \H^k_d$ and $h^k_{d + d^\Lambda} = b_k$.\\
Let $(M,J,\omega,g)$ be an almost K\"ahler $2m$-manifold. Then there is an injection 
    \[
    \H^k_{d + d^c} \longhookrightarrow \H^k_{d + d^\Lambda}.
    \]
    In particular, $h^k_{d + d^c} \le h^k_{d + d^\Lambda}$. Furthermore, if we have $\H^k_{d^c} \subseteq \H^k_{d + d^c}$ for some $k$, then $\H^k_{d + d^c} = \H^k_d$ and $h^k_{d + d^c} = b_k$.
\end{thmx}

In the second part of the paper (sections \ref{sec:surfaces} and \ref{sec:solvmanifolds}), we focus on compact complex surfaces. Recall that for a compact complex manifold of complex dimension $m$, the \emph{Hodge numbers} $h^{p,q}_{\bar \partial}$ are an invariant of the complex structure. They are usually arranged in the so-called \emph{Hodge diamond}
\[
    \begin{matrix}
        & & h^{0,0}_{\bar \partial} & & \\
    & \iddots & \vdots & \ddots & \\
    {h^{m,0}_{\bar \partial}} & \dots & \dots & \dots & {h^{0,m}_{\bar \partial}} \\
    & \ddots & \vdots & \iddots & \\
     & & {h^{m,m}_{\bar \partial}} & &
    \end{matrix}
    \]
The \emph{Bott-Chern numbers} $h^{p,q}_{d+d^c}$ can be arranged in a similar diamond: the \emph{Bott-Chern diamond}. For compact complex surfaces, it turns out that Hodge numbers do not depend on the choice of complex structure, but only on the first Betti number $b_1$ and on the positive and negative self-intersection numbers $b^+$, $b^-$ of the underlying manifold. Since these are topological invariants, we have that Hodge numbers of compact complex surfaces are a consequence of the (oriented) topology of the underlying manifold (cf.\ \cite{BHPV04}). The same is true for Bott-Chern numbers of compact complex surfaces. This fact follows implicitly from the work of Teleman \cite{Tel06}. Alternatively, Stelzig and Wilson \cite[Section 4.1]{SW23} proved the same result using double complexes. Here we explicitly write down the Bott-Chern diamond of compact complex surfaces in terms of topological constants (Theorem \ref{thm:BC:diamond}).

\begin{thmx}
    Let $(M,J)$ be a compact complex surface. If $b_1$ is even, then the Bott-Chern diamond of $(M,J)$ is 
    \[
    \begin{matrix}
        & & 1 & & \\
    & \frac{b_1}{2} & & \frac{b_1}{2} & \\
    \frac{b^+-1}{2} & & b^- +1 & & \frac{b^+-1}{2} \\
    & \frac{b_1}{2} & & \frac{b_1}{2} & \\
     & & 1 & &
    \end{matrix}
    \]
    If $b_1$ is odd, then the Bott-Chern diamond of $(M,J)$ is 
    \[
    \begin{matrix}
        & & 1 & & \\
    & \frac{b_1-1}{2} & & \frac{b_1-1}{2} & \\
    \frac{b^+}{2} & & b^- +1 & & \frac{b^+}{2} \\
    & \frac{b_1+1}{2} & & \frac{b_1+1}{2} & \\
     & & 1 & &
    \end{matrix}
    \]
\end{thmx}

The proof of the theorem is based on results obtained by Teleman \cite{Tel06} on the pseudo-effective cone of non-K\"ahler complex surfaces (see also \cite{ADT16, ATV19}). Note that Aeppli numbers are completely determined by Bott-Chern numbers. From the explicit diamond it is evident that both are topological and independent of the choice of complex structure (\cite[Section 4.1]{SW23} or Corollary \ref{cor:BC:indep}).\\
As a consequence of the theory developed for general almost Hermitian manifolds and for complex surfaces, we are able to determine the numbers $h^k_{d+d^c}$, $h^k_{d+d^\Lambda}$ of $4$-dimensional solvmanifolds endowed with an invariant compatible triple $(J,\omega,g)$, where $J$ is a complex structure (Theorem \ref{thm:solvmanifolds} and tables \ref{tab:abc} and \ref{tab:def}).

\begin{thmx}
    Let $M$ be a $4$-dimensional solvmanifold endowed with an invariant compatible triple $(J,\omega,g)$. Then the numbers $h^k_{d+d^c}$, $h^k_{d + d^\Lambda}$ are independent of the choice of invariant compatible triple.
\end{thmx}

Some of the results valid for solvmanifolds hold more in general for compact quotients of Lie groups by a lattice. We conclude with an application of the theory to the Hopf manifold (Example \ref{ex:hopf}).
\vspace{.2cm}

Finally, we remind that recently Cirici and Wilson \cite{CW21}, Coelho, Placini and Stelzig \cite{CPS22} and the authors \cite{ST21, ST23b}, developed cohomological aspects of the theory which have a strong connection with the theory of harmonic forms presented in this paper.

\section{Preliminaries}\label{sec:preliminaries}
We use the word \emph{manifold} to indicate a compact, connected, smooth manifold with no boundary. We abbreviate $dx^j \wedge dx^k$ to $dx^{jk}$ and $\overline{\phi^j}$ to $\phi^{\bar j}$.
\vspace{.2cm}

An \emph{almost Hermitian $2m$-manifold} $(M,J,\omega,g)$ is a smooth $2m$-manifold $M$ endowed with the following structures:
\begin{itemize}
    \item an \emph{almost complex structure} $J$, i.e., $J \in \End(TM)$, $J^2 = - \Id$;

    \item an \emph{almost symplectic structure} $\omega$, i.e., a non-degenerate, real $2$-form;

    \item a Riemannian metric $g$.
\end{itemize}
These structures must satisfy the compatibility condition
\begin{equation}\label{eq:comp}
    \omega( \cdot \, , J \cdot) = g( \cdot \, , \cdot).
\end{equation}
We say that $(J,\omega, g)$ is a \emph{compatible triple}.\\
Denote by $A^k$ the space of smooth complex $k$-forms on $M$. The existence of a compatible triple induces three operators on $A^k$:
\begin{itemize}
    \item $J$ acts on $k$-forms by duality, inducing a map $J \colon A^k \rightarrow A^k$ and a bigrading decomposition
    \[
    A^k = \bigoplus_{p+q=k} A^{p,q};
    \]
    \item $g$ induces the $\C$-linear Hodge $*$ operator 
        \[
        * \colon A^{p,q} \longrightarrow A^{m-q,m-p};
        \]
        \item $\omega$ induces the symplectic $*_s$ operator $*_s \colon A^{k} \longrightarrow A^{2m-k}$, defined by the relation
        \[
        \alpha \wedge \overline{*_s \, \beta} = \omega( \alpha, \beta) \, \frac{\omega^m}{m!},
        \]
        for any $\alpha, \beta \in A^k$.
\end{itemize}
The compatibility condition \eqref{eq:comp} allows to express $*$ as
\begin{equation}\label{eq:symp:star}
    * = J *_s = *_s J.
\end{equation}
An almost Hermitian manifold can be thought as an almost complex manifold $(M,J)$ together with a fixed $J$-compatible metric $g$. The associated fundamental form $\omega$ is an almost symplectic structure on $M$. On the other side, we can also think of it as an almost symplectic manifold $(M,\omega)$ together with a metric $g$ such that \eqref{eq:comp} holds for some almost complex structure $J$. An almost Hermitian manifold $(M, J , \omega, g)$ is said to be  
\begin{itemize}
    \item \emph{almost K\"ahler} if $d \omega =0$;

    \item \emph{Hermitian} if the Nijenhuis tensor $N_J$ vanishes, i.e., if $J$ is a complex structure;

    \item \emph{K\"ahler} if both $d \omega=0$ and $N_J=0$.
\end{itemize}
On almost symplectic manifolds (thus on almost Hermitian manifolds) it is defined the \emph{Lefschetz operator} $L \colon A^k \rightarrow A^{k+2}$, $\alpha \mapsto \omega \wedge \alpha$, and its \emph{dual Lefschetz operator} $\Lambda \colon A^k \rightarrow A^{k-2}$, $\alpha \mapsto \iota_\omega \alpha$, where $\iota_\omega$ denotes contraction by $\omega$. Since $\omega$ is non-degenerate, for each $k = 1, \dots, m$, powers of the Lefschetz operator give an isomorphism
\[
L^k \colon A^{m-k} \longrightarrow A^{m+k}.
\]
A $k$-form $\alpha$, with $k \le m$, is called \emph{primitive} if $\Lambda \alpha =0$, or, equivalently, if $L^{m-k+1} \alpha =0$. Every $k$-form admits a \emph{Lefschetz decomposition} \cite[Th\'eor\`eme 3, p.\ 26]{Wei58}
\[
\alpha = \sum_{j \ge \max \{ k-m,0 \}}^{} L^j P^{k -2j},
\]
where $P^{k-2j}$ are primitive $(k-2j)$-forms. Given a compatible triple $(J, \omega, g)$, one can compute $*$ on each summand of the Lefschetz decomposition thanks to the formula \cite[Th\'eor\`eme 2, p.\ 23]{Wei58}
\begin{equation}\label{eq:weil}
* L^r P^k = (-1)^{\frac{k (k+1)}{2}} \frac{r!}{(m-k-r)!} L^{m-k-r} JP^k.
\end{equation}

\section{Almost complex and almost symplectic Laplacians on almost Hermitian manifolds}\label{sec:laplacians}

In this section we introduce generalizations of Bott-Chern and Aeppli Laplacians for almost Hermitian manifolds built using almost symplectic differential operators. We also give the definitions of the spaces of $(d+d^\Lambda)$-harmonic forms studied in sections \ref{sec:ddlambda} and \ref{sec:relations}, together with their basic properties. 
\vspace{.2cm}

We begin by reviewing the almost complex point of view, where the main object of study is the space $\H^k_{d + d^c}$ of $(d + d^c)$-harmonic forms. For more details, we refer the reader to \cite{ST23b}, where $(d + d^c)$-harmonic forms have been introduced for the first time on almost Hermitian manifolds.\\
Let $(M,J,\omega,g)$ be an almost Hermitian $2m$-manifold. Consider the operator $d^c := J^{-1} d J$. Then $d^2=0$ and $(d^c)^2 =0$. Observe that $dd^c + d^c d =0$ if and only if $J$ is integrable. Bott-Chern and Aeppli Laplacians for complex manifolds can be generalized to four Laplacians \cite[Definition 4.1]{ST23b}:
{\small
\begin{align*}
    &\Delta_{d + d^c} := d d^c ( d d^c)^* + ( d d^c)^* d d^c + d^* d^c ( d^* d^c )^* + ( d^* d^c )^* d^* d^c + d^* d + (d^c)^* d^c\\
    &\Delta_{d^c + d} := d^c d ( d^c d)^* + ( d^c d)^* d^c d + (d^c)^* d ( (d^c)^* d )^* + ( (d^c)^* d )^* (d^c)^* d + d^* d + (d^c)^* d^c\\
    &\Delta_{d d^c} := d d^c ( d d^c)^* + ( d d^c)^* d d^c +  d (d^c)^* ( d (d^c)^* )^* + ( d (d^c)^* )^* d (d^c)^* + d d^* + d^c (d^c)^* \\
    &\Delta_{d^c d} := d^c d ( d^c d)^* + ( d^c d)^* d^c d + d^c d^* ( d^c d^* )^* + ( d^c d^* )^* d^c d^* + d d^* + d^c (d^c)^*.
\end{align*}
}

If  $P \in \{ d+ d^c, \, d^c +d, \, d d^c, \, d^c d \}$, the \emph{space of $P$-harmonic forms} is
\[
\H^k_P := A^k \cap \ker \Delta_P.
\]
Explicitly, one has that
\begin{align*}
    &\H^k_{d+d^c} = \{ \alpha \in A^k : d \alpha =0, \, d^c \alpha =0, \, (d d^c)^* \alpha =0 \}, \\
    &\H^k_{d^c +d} = \{ \alpha \in A^k : d \alpha =0, \, d^c \alpha =0, \, (d^c d)^* \alpha =0 \}, \\
    &\H^k_{d d^c} = \{ \alpha \in A^k : d^* \alpha =0, \, (d^c)^* \alpha =0, \, d d^c \alpha =0 \}, \\
    &\H^k_{d^c d} = \{ \alpha \in A^k : d^* \alpha =0, \, (d^c)^* \alpha =0, \, d^c d \alpha =0 \}.
\end{align*}
By Hodge theory, the spaces of $P$-harmonic forms are finite dimensional vector spaces over $\C$ \cite[Proposition 4.3]{ST23b} and they are isomorphic to each other according to the following diagram \cite[Proposition 4.5]{ST23b}:
\begin{center}
        \begin{tikzcd}
        \mathcal{H}^k_{d + d^c} \arrow[rr, "\sim", "J"'] \arrow[dd, "\wr", "*"'] &  & \mathcal{H}^k_{d^c + d} \arrow[dd, "\wr", "*"'] \\
                                              &  &                                    \\
\mathcal{H}^{2m-k}_{d d^c} \arrow[rr, "\sim", "J"']         &  & \mathcal{H}^{2m-k}_{d^c d}        
        \end{tikzcd}
        \end{center}
Denoting by $h^k_P$ the complex dimension of $\H^k_P$, one has that
\[
h^k_{d+d^c} = h^k_{d^c + d} = h^{2m-k}_{d d^c} = h^{2m-k}_{d^c d},
\]
and the numbers $h^k_{d + d^c}$ are almost Hermitian invariants \cite[Corollary 4.6]{ST23b}.\\
Let us briefly compare the spaces $\H^k_P$ with the spaces of Bott-Chern and Aeppli harmonic forms on complex manifolds. When $J$ is integrable, we have that $\Delta_{d+d^c} = \Delta_{d^c +d}$ is the Bott-Chern Laplacian, while $\Delta_{d d^c} = \Delta_{d^c d}$ is the Aeppli Laplacian \cite{Sch07}. The spaces of Bott-Chern and Aeppli harmonic forms are defined as
\[
\H^{p,q}_{BC} = A^{p,q} \cap \ker \Delta_{d+d^c} \quad \text{and} \quad \H^{p,q}_{A} = A^{p,q} \cap \ker \Delta_{d d^c}.  
\]
Their dimensions $h^{p,q}_{BC}$, $h^{p,q}_A$ are an invariant of the complex structure, due to an isomorphism with Bott-Chern and Aeppli cohomologies. These numbers are called \emph{Bott-Chern and Aeppli numbers}. In the almost complex case, the operators $dd^c$, $d^c d$ do not preserve the bigrading of forms, thus one has to consider harmonic $k$-forms instead of $(p,q)$-forms (cf.\ \cite[Remark 3.10]{ST23b}).\\
In the next lemma, we show what is the relation between the spaces $\H^k_{d+d^c}$, $\H^k_{d d^c}$ and $\H^{p,q}_{BC}$, $\H^{p,q}_{A}$ in the \emph{integrable} case.

\begin{lemma}\label{lemma:ddc:BC}
    Let $(M,J)$ be a complex manifold. Then
    \[
    \H^k_{d+d^c} = \bigoplus_{p+q=k} \H^{p,q}_{BC} \quad \text{and} \quad \H^k_{d d^c} = \bigoplus_{p+q=k} \H^{p,q}_{A}.
    \]
    In particular, $h^k_{d+d^c} = \sum\limits_{p+q=k} h^{p,q}_{BC}$ and $h^k_{d d^c} = \sum\limits_{p+q=k} h^{p,q}_{A}$.
\end{lemma}
\begin{proof}
    We prove the claim for Bott-Chern harmonic forms. The claim for Aeppli harmonic forms follows by duality. The inclusion $\bigoplus_{p+q=k} \H^{p,q}_{BC} \subseteq \H^k_{d+d^c}$ is immediate since $A^{p,q} \cap \ker \Delta_{d+d^c} \subseteq A^{p+q} \cap \ker \Delta_{d+d^c}$. For the opposite inclusion, let $\alpha \in \H^k_{d+d^c}$ and let 
    \[
    \alpha = \sum_{p+q=k} \alpha^{p,q}
    \]
    be its bidegree decomposition. Since $\alpha$ is both $d$-closed and $d^c$-closed, we have that the forms
    \[
    \alpha^{ev} = \sum_{p \text{ even}} \alpha^{p,q} \quad \text{and} \quad \alpha^{od} = \sum_{p \text{ odd}} \alpha^{p,q}
    \]
    are both $d$-closed (cf.\ \cite[Section 3.3]{ST23b}). Consider the equation
    \[
    0 = d \alpha^{ev} = \bar \partial \alpha^{0,k} + \partial \alpha^{0,k} + \bar \partial \alpha^{2,k-2} + \dots
    \]
    and separate the bidegree of the terms. Since the operators $\partial$, $\bar \partial$ have bidegree $(1,0)$, $(0,1)$ respectively, and since two summands of $\alpha^{ev}$ differ in the first bidegree $p$ by at least $2$, all the terms $\partial \alpha^{p,q}$, $\bar \partial \alpha^{p,q}$, for $p$ even, have different bidegree and for all $(p,q)$ with $p$ even, it must be $\partial \alpha^{p,q} = \bar \partial \alpha^{p,q} =0$. With a similar reasoning, we have that $d \alpha^{p,q} =0$ for all $(p,q)$, showing that $\alpha^{p,q} \in A^{p,q} \cap \ker d \cap \ker d^c$.\\
    The last condition to consider is $(dd^c)^* \alpha=0$. In the complex case, $dd^c$ has bidegree $(1,1)$, thus each summand $\alpha^{p,q}$ is $(dd^c)^*$-closed and $\alpha^{p,q} \in \H^{p,q}_{BC}$, proving the opposite inclusion.
\end{proof}
\vspace{.2cm}

Next, we consider the almost symplectic point of view, studying the space $\H^k_{d + d^\Lambda}$ of $(d + d^\Lambda)$-harmonic forms. This space is introduced here for the first time on almost symplectic manifolds. For the symplectic case, see \cite{TY12a, TY12b, TY14}.\\
Let $d^\Lambda$ be the symplectic co-differential (cf.\ \cite{BT01, Bri88, TY12a}), defined on $k$-forms by
\begin{equation}\label{eq:dlambda}
    d^\Lambda := (-1)^{k+1} *_s d *_s.
\end{equation}
Since $d^2=0$ and $*_s^2 = \Id$, we immediately have that $(d^\Lambda)^2=0$.\\
If $d \omega =0$, then $d d^\Lambda + d^\Lambda d = 0$, but in general $d d^\Lambda + d^\Lambda d \neq 0$. To the best of our knowledge, it is not known if the converse implication is also true, i.e., if $(d d^\Lambda + d^\Lambda d)=0$ implies $d \omega =0$. We are able to prove such implication when $2m = 4$. Indeed, we have the following

\begin{lemma}\label{lemma:balanced}
    Let $(M, \omega)$ be an almost symplectic $2m$-manifold such that $dd^\Lambda + d^\Lambda d=0$. Then $d \omega^{m-1} =0$. In particular, if $2m=4$, then $d \omega =0$.
\end{lemma}
\begin{proof}
    Computing $dd^\Lambda + d^\Lambda d$ on an arbitrary function $f \in C^\infty(M)$, we have that
    \begin{align*}
    0 &= (d d^\Lambda + d^\Lambda d)f = d^\Lambda d f = *_s \, d *_s df = - \frac{1}{(m-1)!} *_s \, d ( \omega^{m-1} \wedge df ) =\\
     &= - \frac{1}{(m-2)!}*_s \, d \omega^{m-1} \wedge df,
    \end{align*}
    where we used \eqref{eq:weil}. Since $*_s$ is an isomorphism, we have that
    \begin{equation}\label{eq:converse}
        d \omega^{m-1} \wedge df =0
    \end{equation}
    for every $f \in C^\infty(M)$. Fix $x \in M$ and let $\{x_j \}_{j =1}^{2m}$ be coordinate functions in a neighborhood $U$ of $x$. In local coordinates we can write
    \[
    d \omega^{m-1} = \sum\limits_{j=1}^{2m} \omega_{j} \, dx^{1 \dots \hat{j} \dots 2m},
    \]
    where $\hat{j}$ denotes missing indices. Choosing $f$ as a smooth extension of $x_j$ to all $M$, \eqref{eq:converse} implies that $\omega_j =0$ for all $j = 1, \dots, 2m$. Therefore
    \[
    d \omega^{m-1} = 0
    \]
    on $U$. Since $x$ is arbitrary, $d \omega^{m-1}=0$ on $M$. If $2m = 4$, then $d \omega=0$. 
\end{proof}

\begin{remark}
    Metrics whose fundamental form satisfies $d \omega^{m-1} =0$ are known in the literature as \emph{(almost) balanced metrics} \cite{Mic82}, or \emph{semi-K\"ahler} \cite{GH80}. Note that in lemma \ref{lemma:balanced} we proved that the balanced condition is equivalent to asking that $d \omega$ is primitive or that $d$ and $d^\Lambda$ anti-commute \emph{on functions}.
\end{remark}

Given a compatible triple $(J, \omega, g)$, by \eqref{eq:symp:star} $d^\Lambda$ has the following expression in terms of $d^c$ \cite{BT01, TY12a}:
\begin{equation}\label{eq:dcdlambda}
d^\Lambda = (d^c)^*.
\end{equation}

On symplectic manifolds, Tseng and Yau \cite{TY12a, TY12b} introduced the Laplacians
{\small
\[
\Delta_{d + d^\Lambda} := d d^\Lambda ( d d^\Lambda)^* + ( d d^\Lambda)^* d d^\Lambda + d^* d^\Lambda ( d^* d^\Lambda )^* + ( d^* d^\Lambda )^* d^* d^\Lambda + d^* d + (d^\Lambda)^* d^\Lambda
\]
}

and
{\small
\[
\Delta_{d d^\Lambda} := d d^\Lambda ( d d^\Lambda)^* + ( d d^\Lambda)^* d d^\Lambda +  d (d^\Lambda)^* ( d (d^\Lambda)^* )^* + ( d (d^\Lambda)^* )^* d (d^\Lambda)^* + d d^* + d^\Lambda (d^\Lambda)^*, 
\]
}

and studied the corresponding spaces of harmonic forms $\H^k_{d + d^\Lambda}$, $\H^k_{d d^\Lambda}$. It turns out that the dimension of such spaces is independent of the choice of metric since they are isomorphic to the symplectic cohomologies $H^k_{d + d^\Lambda}$, $H^k_{d d^\Lambda}$.

\begin{remark}\label{rem:ty}
    Tseng and Yau originally studied \emph{real} harmonic forms. However, since $\Delta_{d + d^\Lambda}$ and $\Delta_{d d^\Lambda}$ are real operators, there is no difference in studying them on complex forms. The resulting spaces of complex harmonic forms will be the complexification of their real counterpart.
\end{remark}
When $d\omega \neq 0$, the operators $d$ and $d^\Lambda$ do not anti-commute and we must consider four different Laplacians.

\begin{definition}\label{def:laplacians}
    The \emph{$(d + d^\Lambda)$-Laplacian} is 
    \begin{equation}\label{def:d+dlambda:laplacian}
    \Delta_{d + d^\Lambda} := d d^\Lambda ( d d^\Lambda)^* + ( d d^\Lambda)^* d d^\Lambda + d^* d^\Lambda ( d^* d^\Lambda )^* + ( d^* d^\Lambda )^* d^* d^\Lambda + d^* d + (d^\Lambda)^* d^\Lambda.
    \end{equation}

    The \emph{$(d^\Lambda + d)$-Laplacian} is 
    \begin{equation}\label{def:dlambda+d:laplacian}
    \Delta_{d^\Lambda + d} := d^\Lambda d ( d^\Lambda d)^* + ( d^\Lambda d)^* d^\Lambda d + (d^\Lambda)^* d ( (d^\Lambda)^* d )^* + ( (d^\Lambda)^* d )^* (d^\Lambda)^* d + d^* d + (d^\Lambda)^* d^\Lambda.
    \end{equation}

    The \emph{$d d^\Lambda$-Laplacian} is 
    \begin{equation}\label{def:ddlambda:laplacian}
    \Delta_{d d^\Lambda} := d d^\Lambda ( d d^\Lambda)^* + ( d d^\Lambda)^* d d^\Lambda +  d (d^\Lambda)^* ( d (d^\Lambda)^* )^* + ( d (d^\Lambda)^* )^* d (d^\Lambda)^* + d d^* + d^\Lambda (d^\Lambda)^* .
    \end{equation}
    
    The \emph{$d^\Lambda d$-Laplacian} is 
    \begin{equation}\label{def:dlambdad:laplacian}
    \Delta_{d^\Lambda d} := d^\Lambda d ( d^\Lambda d)^* + ( d^\Lambda d)^* d^\Lambda d + d^\Lambda d^* ( d^\Lambda d^* )^* + ( d^\Lambda d^* )^* d^\Lambda d^* + d d^* + d^\Lambda (d^\Lambda)^*.
    \end{equation}
    
\end{definition}

In the symplectic case $d$ and $d^\Lambda$ anti-commute, therefore $\Delta_{d + d^\Lambda} = \Delta_{d^\Lambda + d}$ and $\Delta_{d d^\Lambda} = \Delta_{d^\Lambda d}$, recovering the symplectic Laplacians of Tseng and Yau. The symplectic $*_s$ operator and the almost complex structure $J$ provide the following relations among the Laplacians:
\[
*_s \, \Delta_{d+d^\Lambda} = \Delta_{d^\Lambda + d} \, *_s, \quad *_s \, \Delta_{d d^\Lambda} = \Delta_{d^\Lambda d} \, *_s
\]
and
\[
J \Delta_{d+d^\Lambda} = \Delta_{d d^\Lambda} J, \quad J \Delta_{d^\Lambda + d } = \Delta_{d^\Lambda d} J .
\]

\begin{definition}\label{def:P:harmonic}
    Let $P \in \{ d+ d^\Lambda, \, d^\Lambda +d, \, d d^\Lambda, \, d^\Lambda d \}$. A $k$-form $\alpha \in A^k$ is said to be \emph{$P$-harmonic} if $\Delta_P (\alpha) =0$.\\
    We denote the space of $P$-harmonic $k$-forms by $\H^k_P$.
\end{definition}

If $\alpha$ is $P$-harmonic, using the equation $\langle \Delta_P (\alpha), \alpha \rangle =0 $, one can explicitly write the spaces of $P$-harmonic forms. More precisely, we have that
\begin{align*}
    &\H^k_{d+d^\Lambda} = \{ \alpha \in A^k : d \alpha =0, \, d^\Lambda \alpha =0, \, (d d^\Lambda)^* \alpha =0 \}, \\
    &\H^k_{d^\Lambda+d} = \{ \alpha \in A^k : d \alpha =0, \, d^\Lambda \alpha =0, \, (d^\Lambda d)^* \alpha =0 \}, \\
    &\H^k_{d d^\Lambda} = \{ \alpha \in A^k : d^* \alpha =0, \, (d^\Lambda)^* \alpha =0, \, d d^\Lambda \alpha =0 \}, \\
    &\H^k_{d^\Lambda d} = \{ \alpha \in A^k : d^* \alpha =0, \, (d^\Lambda)^* \alpha =0, \, d^\Lambda d \alpha =0 \}.
\end{align*}
We proceed to establish the basic properties of $P$-harmonic forms. First of all, the spaces of harmonic forms are finite-dimensional vector spaces over $\C$.

\begin{proposition}\label{prop:ellipticity}
    Let $P \in \{ d+ d^\Lambda, \, d^\Lambda +d, \, d d^\Lambda, \, d^\Lambda d \}$. Then $\Delta_P$ is a $4^\text{th}$-order, self-adjoint, elliptic operator. There is a decomposition
    \[
    A^k = \H^k_P \overset{\perp}{\oplus} \Ima \Delta_P.
    \]
    Furthermore, $\H^k_P$ is a finite-dimensional vector space over $\C$.
\end{proposition}

\begin{proof}
    The proof closely follows that of \cite[Proposition 4.3]{ST23b}, replacing $d^c$ by $(d^\Lambda)^*$ and $(d^c)^*$ by $d^\Lambda$, thanks to \eqref{eq:dcdlambda}. For the sake of completeness, we write down the core of the proof for $P = d + d^\Lambda$. Denote by $\cong$ the equality up to lower order terms. We have that:
    \begin{itemize}
        \item $d (d^\Lambda)^* + (d^\Lambda)^* d = d d^c + d^c d$ has order one, so that $d (d^\Lambda)^* \cong - (d^\Lambda)^* d$;

        \item $d d^\Lambda +d^\Lambda d$ has order one by the K\"ahler identities for almost Hermitian manifolds \cite{CW20a, FH22}, so that $d d^\Lambda \cong - d^\Lambda d$;

        \item $\Delta_d$ and $\Delta_{d^c}$ are elliptic.
    \end{itemize}
    Thus, we can conclude that 
    \begin{align*}
        \Delta_{d + d^\Lambda} &\cong d d^\Lambda ( d d^\Lambda)^* + ( d d^\Lambda)^* d d^\Lambda + d^* d^\Lambda ( d^* d^\Lambda )^* + ( d^* d^\Lambda )^* d^* d^\Lambda \cong\\
        &\cong d d^* d^\Lambda (d^\Lambda)^* + d^* d ( d^\Lambda)^* d^\Lambda + d^* dd^\Lambda (d^\Lambda )^* + d d^* (d^\Lambda )^* d^\Lambda = \\
        &= \Delta_d \Delta_{d^\Lambda} \cong (\Delta_d)^2,
    \end{align*}
    that is elliptic. The orthogonal direct sum decomposition between image and kernel of $\Delta_P$ and the finite-dimensionality of the kernel follow from the theory of self-adjoint, elliptic operators on compact manifolds.
\end{proof}

\begin{definition}
    Since $\H^k_P$ is finite-dimensional, we set
    \[
    h^k_P := \dim_\C \H^k_P.
    \]
\end{definition} 

The numbers $h^k_P$ are invariants of the almost Hermitian structure. If there is no ambiguity, we will omit dependence of $h^k_P$ on $(J, \omega, g)$.

\begin{proposition}\label{prop:harmonic:iso}
On an almost Hermitian $2m$-manifold $(M,J,\omega, g)$, there is a commutative diagram
        
        \begin{center}
        \begin{tikzcd}
        \mathcal{H}^k_{d + d^\Lambda} \arrow[rr, "\sim", "*_s"'] \arrow[dd, "\wr", "J"'] &  & \mathcal{H}^{2m-k}_{d^\Lambda + d} \arrow[dd, "\wr", "J"'] \\
                                              &  &                                    \\
\mathcal{H}^{k}_{d d^\Lambda} \arrow[rr, "\sim", "*_s"']         &  & \mathcal{H}^{2m-k}_{d^\Lambda d}        
        \end{tikzcd}
        \end{center}
\end{proposition}
\begin{proof}
    The isomorphism between $ \mathcal{H}^k_{d + d^\Lambda}$ and $\mathcal{H}^{2m-k}_{ d^\Lambda + d}$ follows since $\alpha$ is $d$-closed, $d^\Lambda$-closed and $(d d^\Lambda)^*$-closed if and only if $*_s \, \alpha$ is $d^\Lambda$-closed, $d$-closed and $(d^\Lambda d )^*$-closed. The other isomorphisms follow similarly.
\end{proof}

\begin{cor}
Let $(M,J,\omega, g)$ be an almost Hermitian $2m$-manifold. Then
\[
h^k_{d+d^\Lambda} = h^{2m-k}_{d^\Lambda + d} = h^{k}_{d d^\Lambda} = h^{2m-k}_{d^\Lambda d}.
\]
\end{cor}

Note that when $d\omega=0$, the spaces $\H^k_{d + d^\Lambda}$, $\H^k_{d d^\Lambda}$ are isomorphic to the symplectic cohomology groups $H^k_{d + d^\Lambda}$, $H^k_{d d^\Lambda}$. Therefore, their dimensions do not depend on the choice of metric and the numbers $h^k_{d + d^\Lambda}$, $h^k_{d d^\Lambda}$ are symplectic invariants.
\vspace{.2cm}

We conclude this section with a result comparing the numbers $h^k_P$ of almost Hermitian manifolds, for $P \in \{d+ d^c, d^c +d, d d^c, d^c d,d+ d^\Lambda, d^\Lambda +d, d d^\Lambda, d^\Lambda d \}$. Its proof follows closely that of \cite[Theorem 4.4]{ST23b}, therefore it is omitted.

\begin{theorem}
    Let $(M,J,\omega,g)$ and $(M',J',\omega',g')$ be two almost Hermitian manifolds of the same dimension and let $f \colon M \rightarrow M'$ be a surjective smooth map preserving the compatible triple, i.e., satisfying
    \begin{equation}\label{eq:preserving}
    df \circ J = J' \circ df, \quad f^* \omega' = \omega, \quad f^* g' = g.
    \end{equation}
    Then we have that
    \[
    h^k_P (J', \omega', g') \le h^k_P (J, \omega, g).
    \]
    Moreover, if $f$ is also a diffeomorphism then
    \[
    h^k_P (J', \omega', g') = h^k_P (J, \omega, g).
    \]
\end{theorem}

\begin{remark}
    If any two conditions among those in \eqref{eq:preserving} hold, then also the remaining one is satisfied.
\end{remark}

\section{The spaces of \texorpdfstring{$(d+d^c)$}{}-harmonic forms}\label{sec:ddc}

In this section we study the spaces $\H^k_{d + d^c}$ and the numbers $h^k_{d + d^c}$.
\vspace{.2cm}

Let $(M,J,\omega,g)$ be an almost Hermitian $2m$-manifold. We briefly recall the properties of $\H^k_{d+d^c}$ proved in \cite{ST23b}. 
\begin{proposition}[Properties of $\H^k_{d + d^c}$]\label{prop:basic:dlambda}\hfill 

    \begin{itemize}
        \item [(i)] $\H^0_{d + d^c} \cong \H^{2m}_{d + d^c} \cong \C$.

        \item [(ii)] The number $h^1_{d + d^c}$ is metric independent.

        \item [(iii)] If $d\omega =0$ and $2m = 4$, then
        \[
        \H^2_{d + d^c} = \C \langle \omega \rangle \oplus \H^-_g \oplus \H^{(2,0)(0,2)}_J.
        \]
        In particular, $h^2_{d + d^c} = b^- + 1+ h^-_J$ and it does not depend on the choice of the almost K\"ahler metric.
    \end{itemize}
\end{proposition}

We begin by proving a topological upper bound on $h^1_{d + d^c}$.

\begin{lemma}\label{lemma:h1bound}
    Let $(M,J)$ be an almost complex $2m$-manifold. Then
    \[
    h^1_{d + d^c} \le b_1.
    \]
\end{lemma}
\begin{proof}
    We prove that $\H^1_{d + d^c}$ injects into de Rham cohomology. From the definition, it follows that
    \[
    \H^1_{d+d^c} = ( A^{1,0} \cap \ker d) \cup ( A^{0,1} \cap \ker d).
    \]
    Let $\alpha \in \H^1_{d + d^c}$. Then $\alpha = \alpha^{1,0} + \alpha^{0,1}$, with $\alpha^{1,0}$, $\alpha^{0,1}$ both $d$-closed, and $\alpha$ defines a de Rham class. To prove injectivity, suppose that
    \[
    \alpha^{1,0} + \alpha^{0,1} = df
    \]
    for some $f \in C^\infty (M)$. By bidegree, we have that $\alpha^{1,0} = \partial f$ and $\alpha^{0,1} = \bar \partial f$. Since $d \alpha^{1,0} =0$, we have that $\partial^2 f + \bar \partial \partial f + \bar \mu \partial f =0$. In particular, by bidegree reasons $\bar \partial \partial f =0$, which implies that $f$ is constant, thus $\alpha^{1,0} = \partial f =0$ and $\alpha^{0,1} = \bar \partial f =0$.
\end{proof}

\begin{cor}\label{cor:h1:b1}
    Let $(M,J)$ be an almost complex $2m$-manifold such that $b_1 \in \{ 0, 1 \}$. Then $h^1_{d + d^c}=0$.
\end{cor}
\begin{proof}
    By lemma \ref{lemma:h1bound}, we have $h^1_{d + d^c} \le b_1 \le 1$. Moreover, $h^1_{d + d^c}$ is even because it is invariant under conjugation.
\end{proof}
We invite the reader to compare the results of lemma \ref{lemma:h1bound} and corollary \ref{cor:h1:b1} with \cite[Corollary 4.6]{CW20b} and \cite[Lemma 4.2]{HPT23}.

We now establish a decomposition of $\H^2_{d + d^c}$ valid on almost Hermitian $4$-manifolds. Let $(M,J,\omega,g)$ be an almost Hermitian $4$-manifold. Consider the harmonic decomposition of $\omega$ with respect to the Hodge Laplacian $\Delta_d$, i.e.,
\begin{equation}\label{eq:harmonic:omega}
\omega = h(\omega) + d \eta + d^* \mu,
\end{equation}
where $h(\omega)$ is $d$-harmonic, $\eta \in A^1$, and $\mu \in A^3$. Define the $2$-form $\gamma_0$ as
\begin{equation}\label{eq:gamma0}
    \gamma_0 := -d*\mu -d^* \mu.
\end{equation}
Note that 
\begin{equation}\label{eq:gamma0:*}
    * \gamma_0 = - *d*\mu + *^2 d *\mu = d^* \mu + d* \mu = - \gamma_0,
\end{equation}
hence $\gamma_0$ is anti-self-dual. Since anti-self-dual forms have necessarily bidegree $(1,1)$, we also have that
\begin{equation}\label{eq:gamma0:J}
    J \gamma_0 = \gamma_0.
\end{equation}
Denote by $\H^-_g$ the space of anti-self-dual, $d$-harmonic forms and by $\H^{(2,0)(0,2)}_J$ the space of $J$-anti-invariant, $d$-harmonic forms.\\
We prove the following

\begin{theorem}\label{thm:ddc:dec}
Let $(M,J,\tilde{\omega},\tilde{g})$ be an almost Hermitian $4$-manifold and let  $(\omega,g)$ be a Gauduchon metric in the same conformal class of $(\tilde{\omega}, \tilde{g})$. Then
\[
\H^2_{d + d^c} = \C \langle \omega + \gamma_0 \rangle \oplus \H^-_g \oplus \H^{(2,0)(0,2)}_J.
\]
In particular, $h^2_{d+d^c} = b^- +1 + h^-_J$ and it is metric independent.
\end{theorem}
\begin{proof}
Let $g$ be a metric in the same conformal class of $\Tilde{g}$. If $\alpha \in A^2$, then $ *_g \alpha =*_{\tilde{g}} \alpha$. As a consequence, the space $\H^2_{d + d^c}$ is invariant under conformal changes of metric. In each conformal class of metric there is always a Gauduchon metric \cite{Gau77}, i.e., a metric for which $d d^c \omega = d^c d \omega = 0$. Therefore, we assume that $g$ is a Gauduchon metric in the same conformal class of $\Tilde{g}$.\\
We first prove the inclusion $\C \langle \omega + \gamma_0 \rangle \oplus \H^-_g \oplus \H^{(2,0)(0,2)}_J \subseteq \H^2_{d + d^c}$.\\
We have that $\omega + \gamma_0$ is $d$-closed since 
\[
d(\omega + \gamma_0) = d ( h(\omega) + d\eta - d* \mu) = 0,
\]
it is $d^c$-closed by \eqref{eq:gamma0:J}, and it is $(d d^c)^*$-closed by \eqref{eq:gamma0:*} and because the metric is Gauduchon.
Forms in $\H^-_g \oplus \H^{(2,0)(0,2)}_J$ are necessarily $(d+d^c)$-harmonic since they are $d$-harmonic and have bidegree either $(1,1)$ or $(2,0)+(0,2)$.\\
We now prove the opposite inclusion. Let $\alpha \in \H^2_{d + d^c}$. Write $\alpha$ using the Lefschetz decomposition and the bidegree as
\[
\alpha = f \omega + \gamma^{1,1} + \gamma^{(2,0)(0,2)},
\]
with $\gamma^{1,1}$, $\gamma^{(2,0)(0,2)}$ primitive forms and $f \in C^\infty(M)$. Since $d\alpha=0$ and $d^c \alpha =0$, we have that
\begin{equation}\label{eq:domega:gamma}
d(f \omega + \gamma^{1,1})=0
\end{equation}
and
\begin{equation}\label{eq:dgamma2}
d\gamma^{(2,0)(0,2)} =0.
\end{equation}
The form $\gamma^{(2,0)(0,2)}$ is $d$-closed and primitive, thus $d$-harmonic, which implies that $\gamma^{(2,0)(0,2)} \in \H^{(2,0)(0,2)}_J$. On the other side, we have that 
\begin{align*}
    0 = (d d^c)^* \alpha &= - * d^c d * (f \omega + \gamma ^{1,1} + \gamma^{(2,0)(0,2)})=\\
    &= - * d^c d (f \omega - \gamma ^{1,1} + \gamma^{(2,0)(0,2)})=\\
    &= - 2 * d^c d (f \omega),
\end{align*}
where in the last equality we used \eqref{eq:domega:gamma} and \eqref{eq:dgamma2}. Since the metric is Gauduchon, we have that 
\[
0 = d^c d (f \omega) = d^c d f \wedge \omega - df \wedge d^c \omega + d^c f \wedge d\omega.
\] 
Consider the real operator $P \colon C^\infty(M) \rightarrow C^\infty(M)$ given by 
\[
P(f) = *(d^c d f \wedge \omega - df \wedge d^c \omega + d^c f \wedge d\omega).
\]
By the same argument of \cite[Proof of Theorem 4.3]{PT22}, $P$ is strongly elliptic, thus $f$ must be constant. Consider the form
\[
\beta^{1,1} := \gamma^{1,1} - f \gamma_0.
\]
The form $\beta^{1,1}$ is anti-self-dual by \eqref{eq:gamma0:*}, and it is also $d$-closed by \eqref{eq:harmonic:omega}, \eqref{eq:gamma0} and \eqref{eq:domega:gamma}. Therefore it is also $d$-harmonic and $\beta^{1,1} \in \H^-_g$. The claim follows writing $\alpha$ as
\[
\alpha = f (\omega+\gamma_0) + \beta^{1,1} + \gamma^{(2,0)(0,2)},
\]
with $f$ constant, $\beta^{1,1} \in \H^-_g$ and $\gamma^{(2,0)(0,2)} \in \H^{(2,0)(0,2)}_J$.
\end{proof}

\begin{remark}
    Applying theorem \ref{thm:ddc:dec} to almost K\"ahler $4$-manifolds, we recover \cite[Theorem 5.1]{ST23b}, since when $d \omega=0$, we have $\gamma_0 =0$. Our result should be compared with \cite[Theorem 4.2]{Hol22}.
\end{remark}

We conclude the section with a result that in certain cases allows to explicitly compute $h^k_{d + d^c}$ on almost K\"ahler manifolds of any dimension. 

\begin{theorem}\label{thm:ak:betti}
    Let $(M,J, \omega,g)$ be an almost K\"ahler $2m$-manifold. Suppose that we have $\H^k_{d^c} \subseteq \H^k_{d + d^c}$ for some $k$. Then $\H^k_{d + d^c} = \H^k_d$ and $h^k_{d + d^c} = b_k$.
\end{theorem}
\begin{proof}
    Let $\alpha \in \H^k_{d + d^c}$ and let 
    \[
    \alpha = h_{d^c} (\alpha) + d^c \eta + (d^c)^* \gamma
    \]
    be its Hodge decomposition with respect to $\Delta_{d^c}$. Since $\alpha$ is $d^c$-closed, we have $(d^c)^*\gamma =0$. The form
    \[
    d^c \eta = \alpha - h_{d^c} (\alpha)
    \]
    is $(d + d^c)$-harmonic because it is the difference of two harmonic forms by the assumption on $\H^k_{d^c}$. From the equation $(d d^c)^* \alpha =0$, we deduce that
    \[
    (d^c)^* d^* d^c \eta = (d^c)^* d^* (\alpha - h_{d^c}(\alpha) ) = 0.
    \]
    By \eqref{eq:dcdlambda} and $d \omega=0$, we have that
    \[
    0 = (d^c)^* d^* d^c \eta = d^\Lambda d^* (d^\Lambda)^* \eta = - d^\Lambda (d^\Lambda)^* d^* \eta,
    \]
    which implies $\langle d^\Lambda (d^\Lambda)^* d^* \eta, d^* \eta \rangle =0$. Therefore $(d^\Lambda)^* d^* \eta= - d^* (d^\Lambda)^* \eta =0$. Finally, $d \alpha=0$ since $\alpha \in \H^k_{d + d^\Lambda}$ and $d^* \alpha = d^* (d^\Lambda)^* \eta=0$, so that $\alpha \in \H^k_d$. We proved that $\H^k_{d+d^c} \subseteq \H^k_d$. To conclude, observe that
    \[
    \H^k_{d+d^c} \subseteq \H^k_d \cong \H^k_{d^c} \subseteq \H^k_{d+d^c},
    \]
    giving the equality of the spaces.
\end{proof}

\section{The spaces of \texorpdfstring{$(d+d^\Lambda)$}{}-harmonic forms}\label{sec:ddlambda}

In this section we study the spaces $\H^k_{d + d^\Lambda}$ and the numbers $h^k_{d + d^\Lambda}$.
\vspace{.2cm}

We begin by proving some basic results valid in any dimension.\\
Let $(M,J,g,\omega)$ be an almost Hermitian $2m$-manifold.

\begin{proposition}[Properties of $\H^k_{d + d^\Lambda}$]\label{prop:basic} \hfill 
    \begin{itemize}
        \item [(i)] $\H^0_{d + d^\Lambda} \cong \H^{2m}_{d + d^\Lambda} \cong \C$.

        \item [(ii)] There is an inclusion $\H^1_{d + d^\Lambda} \subseteq \H^{1}_{d}$. In particular, $h^1_{d + d^\Lambda} \le b_1$.

        \item [(iii)] If $d\omega =0$, then $\H^1_{d + d^\Lambda} = \H^{1}_{d}$ and $h^1_{d + d^\Lambda} = b_1$.
    \end{itemize}
\end{proposition}
\begin{proof}
    We prove (i). Note that any function $f \in \H^0_{d + d^\Lambda}$ must be $d$-closed and therefore constant. Let $f \, Vol \in \H^{2m}_{d + d^\Lambda}$. Since $ d^\Lambda (f \, Vol) = (d^c)^* (f \, Vol) =0$, we have that
    \[
    0= d J * (f \, Vol) = df,
    \]
    so that $\H^{2m}_{d + d^\Lambda} = \C \langle Vol \rangle \cong \C$. To prove (ii), let $\alpha \in \H^{1}_{d + d^\Lambda}$. Then $\alpha$ is $d$-closed and $(dd^\Lambda)^*$-closed. From the equation $(dd^\Lambda)^* \alpha =0$, we have that
    \[
    0 = d J * d * \alpha  = - d J d^* \alpha = - d d^* \alpha,
    \]
    since $J$ acts trivially on functions. Therefore $d^* \alpha =0$ and $\alpha$ is $d$-harmonic. If $d \omega =0$, the opposite inclusion also holds. Indeed, let $\alpha \in \H^1_d$. Then we immediately have $d \alpha =0$ and $ ( d d^\Lambda)^* \alpha = d^c d^* \alpha =0$. Moreover, since $d \omega =0$, we also have $d^\Lambda \alpha = (d \Lambda - \Lambda d) \alpha =0$ since $\alpha$ is $d$-closed and has degree $1$, proving (iii).
\end{proof}

\begin{remark}
    Point (iii) of proposition \ref{prop:basic} is a well-known fact (cf.\ \cite[Lemma 2.7]{BT01}).
\end{remark}

On $4$-manifolds, the study of $(d + d^\Lambda)$-harmonic $2$-forms is reduced to primitive forms.

\begin{theorem}\label{thm:ddlambda:dec}
    Let $(M,J,\omega,g)$ be an almost Hermitian $4$-manifold. Let $\mathcal{PH}^2_{d + d^\Lambda}$ be the space of $(d + d^\Lambda)$-harmonic primitive $2$-forms. Then
    \[
    \H^2_{d + d^\Lambda} = 
    \begin{cases}
        \C \langle \omega \rangle \oplus 
    \mathcal{PH}^2_{d + d^\Lambda} & \text{if } d \omega = 0,\\
    \mathcal{PH}^2_{d + d^\Lambda} & \text{if } d \omega \neq 0.
    \end{cases}
    \]
\end{theorem}

\begin{proof} The inclusions $\C \langle \omega \rangle \oplus 
    \mathcal{PH}^2_{d + d^\Lambda} \subseteq \H^2_{d + d^\Lambda}$ and $\mathcal{PH}^2_{d + d^\Lambda} \subseteq \H^2_{d + d^\Lambda}$ are immediate. For the opposite inclusions, let $\alpha \in \H^2_{d + d^\Lambda}$ and let
    \[
    \alpha = f \omega + \gamma^{1,1} + \gamma^{(2,0)(0,2)}
    \]
    be its Lefschetz and bidegree decomposition. Since $\alpha$ is $d$-closed and $d^\Lambda$-closed, we have that
    \begin{equation}\label{eq:dfomega}
    0 = d J * \alpha = d( f \omega - \gamma^{1,1} - \gamma^{(2,0)(0,2)}) = d ( 2f \omega - \alpha) = 2 d (f \omega).
    \end{equation}
    Thus $d d^c ( f\omega ) = d^c d (f\omega) =0 $ and with the same argument of the proof of theorem \ref{thm:ddc:dec}, we deduce that $f$ is constant. Since $f \omega$ is $d$-closed, $d^\Lambda$ closed and self-dual, it is $(d + d^\Lambda)$-harmonic. In particular, the form $\gamma^{1,1} + \gamma^{(2,0)(0,2)}$ is primitive and $(d + d^\Lambda)$-harmonic. If $d \omega =0$, we obtain that
    \[
    \H^2_{d + d^\Lambda} = \C \langle \omega \rangle \oplus 
    \mathcal{PH}^2_{d + d^\Lambda},
    \]
    proving the first part of the theorem. For the second part, by \eqref{eq:dfomega} we have that
    \[
    0 = d( f\omega ) = f d\omega,
    \]
    with $f$ constant. In particular, if $d \omega \neq 0$, then $f =0$.
\end{proof}

We conclude the section with the Hermitian counterpart of theorem \ref{thm:ak:betti}, which allows to explicitly compute the numbers $h^k_{d + d^\Lambda}$ when $J$ is integrable.

\begin{theorem}\label{thm:hermitian:betti}
    Let $(M,J, \omega,g)$ be a Hermitian $2m$-manifold. Suppose that we have $\H^k_{d^\Lambda} \subseteq \H^k_{d + d^\Lambda}$ for some $k$. Then $\H^k_{d + d^\Lambda} = \H^k_d$ and $h^k_{d + d^\Lambda} = b_k$.
\end{theorem}
\begin{proof}
    The proof follows closely that of theorem \ref{thm:ak:betti}, replacing $d^c$ by $(d^c)^*$. Let $\alpha \in \H^k_{d + d^\Lambda}$ and let 
    \[
    \alpha = h_{d^\Lambda} (\alpha) + d^\Lambda \eta + (d^\Lambda)^* \gamma
    \]
    be its Hodge decomposition with respect to $\Delta_{d^\Lambda}$. Since $\alpha$ is $d^\Lambda$-closed, we have $(d^\Lambda)^*\gamma =0$. From the equation $(d d^\Lambda)^* \alpha =0$, we deduce that
    \[
    d^c d^* d^\Lambda \eta = d^c d^* ( \alpha - h_{d^\Lambda}(\alpha)) = 0.
    \]
    By \eqref{eq:dcdlambda} and integrability of $J$, we have
    \[
    0 = d^c d^* d^\Lambda \eta = d^c d^* (d^c)^* \eta = - d^c (d^c)^* d^* \eta,
    \]
    which implies $(d^c)^* d^* \eta = - d^* (d^c)^* \eta =0$. Since $d \alpha=0$, $d^* \alpha = d^* (d^c)^* \eta=0$, and
    \[
    \H^k_{d+d^\Lambda} \subseteq \H^k_d \cong \H^k_{d^\Lambda} \subseteq \H^k_{d+d^\Lambda},
    \]
    the theorem is proved.
\end{proof}

\section{Relations between \texorpdfstring{$(d + d^c)$}{}-harmonic forms and \texorpdfstring{$(d + d^\Lambda)$}{}-harmonic forms}\label{sec:relations}

In this section we compare the spaces $\H^k_{d + d^c}$, $\H^k_{d + d^\Lambda}$ and the numbers $h^k_{d + d^c}$, $h^k_{d + d^\Lambda}$ on almost Hermitian manifolds.

The first result is an injection of $(d +d^c)$-harmonic forms into $(d +d^\Lambda)$-harmonic forms valid for almost K\"ahler manifolds.

\begin{theorem}\label{thm:akinclusion}
    Let $(M,J,\omega,g)$ be an almost K\"ahler $2m$-manifold. There is an injection 
    \[
    \H^k_{d + d^c} \longhookrightarrow \H^k_{d + d^\Lambda}.
    \]
    In particular, $h^k_{d + d^c} \le h^k_{d + d^\Lambda}$. 
\end{theorem}
\begin{proof}
    Let $\alpha \in \H^k_{d + d^c}$. Since $\alpha$ is $d$-closed, $d^c$-closed and $(d d^c)^*$-closed, we have that 
    \[
    0 = (d^c d)^* (J \alpha) = d^* d^\Lambda (J \alpha) = d^* (d \Lambda - \Lambda d) (J \alpha) = d^* d (\Lambda J \alpha),
    \]
    which implies $d\Lambda (J\alpha ) = d^\Lambda (J\alpha)=0$. Furthermore, we also have $d( J \alpha)=0$. In particular, $J \alpha$ is both $d$-closed and $d^\Lambda$-closed, thus it defines a symplectic cohomology class $[J \alpha]_{d + d^\Lambda} \in H^k_{d + d^\Lambda}$. Taking the harmonic representative $h_{d + d^\Lambda} ( J \alpha)$ we have a well-defined map $\H^k_{d + d^c} \rightarrow \H^k_{d + d^\Lambda}$. The map is injective because if $ J \alpha = d^\Lambda d \beta$ for some $\beta \in A^k$, then 
    \[
    0 = - J d \alpha = d^c (J \alpha) = ( d^\Lambda)^* d^ \Lambda d \beta,
    \]
    giving $J \alpha = d^\Lambda d \beta =0$.
\end{proof}

The opposite inclusion in general does not hold. For instance, one can endow the Kodaira-Thurston manifold with an almost K\"ahler structure such that $h^1_{d + d^c} = 2$ and $h^1_{d + d^\Lambda} = b_1 = 3$ \cite[Proposition 6.2]{ST23b}. Nevertheless, we can prove the opposite inclusion on $(2m-1)$-forms.

\begin{theorem}\label{thm:2m-1}
    Let $(M,J,g,\omega)$ be an almost Hermitian $2m$-manifold. Then
    \[
    \H^{2m-1}_{d^\Lambda + d} \subseteq \H^{2m-1}_{d^c + d}.
    \]
    If $d \omega =0$, then $\H^{2m-1}_{d + d^\Lambda} = \H^{2m-1}_{d^\Lambda + d} = \H^{2m-1}_{d + d^c}$.
\end{theorem}
\begin{proof}
    Let $\alpha \in \H^{2m-1}_{d^\Lambda + d}$. Since $\alpha$ is $(d^\Lambda d)^*$-closed, we have that
    \[
    0 = d^* J^{-1} d J \alpha = d^* d J \alpha,
    \]
    where we used the fact that $J = \Id$ on top-forms. Hence, we conclude that $d^c\alpha=0$. By the equation $d^\Lambda \alpha=0$, we have that $d J * \alpha = 0$. Thus $J \alpha$ is $(d + d^c)$-harmonic, since
    \[
    (dd^c)^* J \alpha = d^\Lambda d^* J \alpha = 0,
    \]
    and $\alpha$ is $(d^c + d)$-harmonic. This proves the inclusion $\H^{2m-1}_{d^\Lambda + d} \subseteq \H^{2m-1}_{d^c + d}$. If $d \omega =0$, by theorem \ref{thm:akinclusion} we have that
    \[
    \H^{2m-1}_{d^\Lambda +d} \subseteq \H^{2m-1}_{d^c +d} \cong \H^{2m-1}_{d +d^c} \hookrightarrow \H^{2m-1}_{d +d^\Lambda} = \H^{2m-1}_{d^\Lambda +d},
    \]
    giving the equality of the spaces and concluding the proof.
\end{proof}

\begin{remark}
Applying theorem \ref{thm:2m-1} to almost K\"ahler $4$-manifolds, we recover \cite[Theorem 5.4]{ST23b}.
\end{remark}

\begin{cor}\label{cor:2m-1}
    Let $(M,J,\omega,g)$ be an almost Hermitian $2m$-manifold. Then $h^1_{d+d^\Lambda} \le h^{2m-1}_{d+d^c}$. If $d \omega =0$, then $h^{2m-1}_{d+d^c} = b_1$.
\end{cor}

\section{Bott-Chern and Aeppli numbers of compact complex surfaces}\label{sec:surfaces}

The goal of this section is to compute Bott-Chern and Aeppli numbers of compact complex surfaces. As a consequence of the computations, we see that they depend only on the topology of the underlying manifold. More precisely, they do not depend on the choice of complex structure, but only on the numbers $b_1$, $b^+$, $b^-$. This is a result which was already implicitly contained in the work of Teleman \cite{Tel06} (see also \cite{SW23}).
\vspace{.2cm}

Let $(M,J)$ be a compact complex surface (with no boundary). We are interested in the following invariants:
\begin{itemize}
    \item the \emph{Betti numbers} of $M$,
    \[
    b_k = \dim_\C H^k_d (M; \C), \quad k=0,\dots,4.
    \]
    Since $M$ is an oriented closed manifold, the Betti numbers reduce to
    \[
    \begin{array}{ c|c|c|c|c|c } 
    k & 0 & 1 & 2 & 3 & 4 \\ 
    \hline
    b_k & 1 & b_1 & b_2 & b_1 & 1 \\
    \end{array}
    \]

    \item the \emph{Hodge numbers} of $(M,J)$,
    \[
    h^{p,q}_{\bar \partial} = \dim_\C H^{p,q}_{\bar \partial} (M; \C), \quad p,q = 0,1,2,
    \]
    that can be arranged in the so called \emph{Hodge diamond}
    \[
    \begin{matrix}
        & & h^{0,0}_{\bar \partial} & & \\
    & {h^{1,0}_{\bar \partial}} & & {h^{0,1}_{\bar \partial}} & \\
    {h^{2,0}_{\bar \partial}} & & {h^{1,1}_{\bar \partial}} & & {h^{0,2}_{\bar \partial}} \\
    & {h^{2,1}_{\bar \partial}} & & {h^{1,2}_{\bar \partial}} & \\
     & & {h^{2,2}_{\bar \partial}} & &
    \end{matrix}
    \]

    \item the \emph{Bott-Chern and Aeppli numbers} of $(M,J)$
    \begin{alignat*}{2}
    &h^{p,q}_{BC} = \dim_\C H^{p,q}_{BC} (M; \C), \quad &p,q = 0,1,2, \\
    &h^{p,q}_{A} = \dim_\C H^{p,q}_{A} (M; \C), \quad &p,q = 0,1,2.
    \end{alignat*}
    By duality between Bott-Chern and Aeppli cohomology, we have that $h^{p,q}_{BC} = h^{m-p,m-q}_{A}$ for all $p$, $q$, thus knowing Bott-Chern numbers completely determines Aeppli numbers. We arrange Bott-Chern numbers in a \emph{Bott-Chern diamond}
    \[
    \begin{matrix}
        & & h^{0,0}_{BC} & & \\
    & {h^{1,0}_{BC}} & & {h^{0,1}_{BC}} & \\
    {h^{2,0}_{BC}} & & {h^{1,1}_{BC}} & & {h^{0,2}_{BC}} \\
    & {h^{2,1}_{BC}} & & {h^{1,2}_{BC}} & \\
     & & {h^{2,2}_{BC}} & &
    \end{matrix}
    \]
\end{itemize}

It is a well-known fact that while a priori Hodge numbers depend on the choice of $J$, for compact complex surfaces they actually depend only on the first Betti number $b_1$ and on the positive and negative self-intersection numbers $b^+$ and $b^-$, with $b^+ + b^- = b_2$. For the sake of completeness, we give here a precise statement, whose proof follows from \cite[Chapter 4, Theorems 2.7 and 2.14]{BHPV04}. 

\begin{theorem}(\cite{BHPV04})\label{thm:hodge:diamond}
    Let $(M,J)$ be a compact complex surface. If $b_1$ is even, then the Hodge diamond of $(M,J)$ is 
    \[
    \begin{matrix}
        & & 1 & & \\
    & \frac{b_1}{2} & & \frac{b_1}{2} & \\
    \frac{b^+-1}{2} & & b^- +1 & & \frac{b^+-1}{2} \\
    & \frac{b_1}{2} & & \frac{b_1}{2} & \\
     & & 1 & &
    \end{matrix}
    \]
    If $b_1$ is odd, then the Hodge diamond of $(M,J)$ is 
    \[
    \begin{matrix}
        & & 1 & & \\
    & \frac{b_1-1}{2} & & \frac{b_1+1}{2} & \\
    \frac{b^+}{2} & & b^- & & \frac{b^+}{2} \\
    & \frac{b_1+1}{2} & & \frac{b_1-1}{2} & \\
     & & 1 & &
    \end{matrix}
    \]
\end{theorem}

\begin{cor}
    Hodge numbers of compact complex surfaces depend only on the topology of the manifold.
\end{cor}

We proceed now to state and prove a similar result, valid for Bott-Chern numbers. 

\begin{theorem}\label{thm:BC:diamond}
    Let $(M,J)$ be a compact complex surface. If $b_1$ is even, then the Bott-Chern diamond of $(M,J)$ is 
    \[
    \begin{matrix}
        & & 1 & & \\
    & \frac{b_1}{2} & & \frac{b_1}{2} & \\
    \frac{b^+-1}{2} & & b^- +1 & & \frac{b^+-1}{2} \\
    & \frac{b_1}{2} & & \frac{b_1}{2} & \\
     & & 1 & &
    \end{matrix}
    \]
    If $b_1$ is odd, then the Bott-Chern diamond of $(M,J)$ is 
    \[
    \begin{matrix}
        & & 1 & & \\
    & \frac{b_1-1}{2} & & \frac{b_1-1}{2} & \\
    \frac{b^+}{2} & & b^- +1 & & \frac{b^+}{2} \\
    & \frac{b_1+1}{2} & & \frac{b_1+1}{2} & \\
     & & 1 & &
    \end{matrix}
    \]
\end{theorem}
\begin{proof}
    One easily sees that $h^{0,0}_{BC} = h^{2,2}_{BC} = 1$ for every compact complex surface. By \cite{Tel06}, see also \cite{ADT16, ATV19}, on compact complex surfaces we have that
    \begin{equation}\label{eq:delta1}
        h^{1,0}_{BC} + h^{0,1}_{BC} + h^{2,1}_{BC} + h^{1,2}_{BC} = 2 b_1
    \end{equation}
    and
    \begin{equation}\label{eq:delta2}
        h^{2,0}_{BC} + h^{1,1}_{BC} + h^{0,2}_{BC} = 
        \begin{cases}
            b_2 & \text{if $b_1$ is even,}\\
            b_2 +1 & \text{if $b_1$ is odd.}
        \end{cases}
    \end{equation}
    Since Bott-Chern numbers are symmetric in $p$,$q$, we can simplify \eqref{eq:delta1} and \eqref{eq:delta2} to get
    \begin{equation}\label{eq:new:delta1}
        h^{1,0}_{BC} + h^{2,1}_{BC} = b_1
    \end{equation}
    and
    \begin{equation}\label{eq:new:delta2}
        2 h^{2,0}_{BC} + h^{1,1}_{BC} = 
        \begin{cases}
            b_2 & \text{if $b_1$ is even,}\\
            b_2 +1 & \text{if $b_1$ is odd.}
        \end{cases}
    \end{equation}
    If $b_1$ is even, then $(M,J)$ admits a K\"ahler metric \cite{Buc99,Lam99}. Since $h^{p,q}_{BC}$ are independent of the choice of metric, it is enough to compute them for a K\"ahler metric. On K\"ahler manifolds, there is an isomorphism $H^{p,q}_{BC} \cong H^{p,q}_{\bar \partial}$. Therefore the Bott-Chern diamond of compact complex surfaces with even $b_1$ coincides with their Hodge diamond given in theorem \ref{thm:hodge:diamond}.\\
    Suppose now that $b_1$ is odd. For any choice of metric $g$, we have $\H^{1,0}_{BC} = \H^{1,0}_{\bar \partial}$. Indeed, writing explicitly the spaces of harmonic forms, we have that
    \[
    \H^{1,0}_{BC} = A^{1,0} \cap \ker \partial \cap \ker \bar \partial \quad \text{and} \quad \H^{1,0}_{\bar \partial} = A^{1,0} \cap \ker \bar \partial.
    \]
    By \cite[Chapter 4, Lemma 2.1]{BHPV04}, every holomorphic form on a compact complex surface is $d$-closed, thus $A^{1,0} \cap \ker \bar \partial = A^{1,0} \cap \ker \partial \cap \ker \bar \partial$, giving the equality of the two spaces, and allowing to deduce that
    \[
    h^{1,0}_{BC} = h^{0,1}_{BC} = h^{1,0}_{\bar \partial} = \frac{b_1 -1}{2}.
    \]
    By \eqref{eq:new:delta1}, we also obtain
    \[
    h^{2,1}_{BC} = h^{1,2}_{BC} = \frac{b_1 + 1}{2}.
    \]
    The number $h^{1,1}_{BC}$ can be computed either using \cite[Lemma 2.3]{Tel06}, using theorem \ref{thm:ddc:dec} applied to an integrable $J$ together with proposition \ref{lemma:ddc:BC}, or applying \cite[Theorem 4.2]{Hol22}, since in the complex case $d$, $d^c$ or $\partial$, $\bar \partial$ are interchangeable in the definition of Bott-Chern cohomology. It turns out that for Bott-Chern numbers
    \begin{equation}\label{eq:h11}
    h^{1,1}_{BC} = b^- +1.
    \end{equation}
    Finally, by \eqref{eq:new:delta2} and \eqref{eq:h11}, we have
    \[
    h^{2,0}_{BC} = h^{0,2}_{BC} = \frac{b^+}{2},
    \]
    concluding the proof.
\end{proof}

\begin{remark} Note that, for Bott-Chern numbers, one has $h^{1,1}_{BC} = b^- +1$ independently of the parity of $b_1$, in contrast to what happens for the Hodge number $h^{1,1}_{\bar \partial}$ (cf.\ theorem \ref{thm:hodge:diamond}). 
\end{remark}

\begin{cor}[{\cite{Tel06} or \cite{SW23}}]\label{cor:BC:indep}
    Bott-Chern and Aeppli numbers of compact complex surfaces depend only on the topology of the underlying manifold.
\end{cor}

\section{The numbers \texorpdfstring{$h^k_{d+d^c}$, $h^k_{d+d^\Lambda}$}{} of invariant compatible triples}\label{sec:solvmanifolds}

In this section we determine the numbers $h^k_{d+d^c}$, $h^k_{d + d^\Lambda}$ of compact complex surfaces that are diffeomorphic to solvmanifolds endowed with an invariant compatible triple. Along the way, we establish results valid more in general for compact quotients of Lie groups by a lattice.
\vspace{.2cm}

Let $G$ be a Lie group and let $\Gamma \subset G$ be a lattice. Suppose that the quotient $M = \Gamma \backslash G$ is a compact manifold. When $G$ is solvable/nilpotent, $M$ is called a \emph{solvmanifold/nilmanifold}. An almost complex structure on $M$ is \emph{invariant} if it is induced by a complex structure on $\g$, the Lie algebra of $G$. Similarly, a metric and an almost symplectic structure on $M$ are \emph{invariant} if they are induced by an inner product on $\g$ and a non-degenerate element of $\bigwedge^2 \g^*$, respectively. An \emph{invariant compatible triple} on $M$ is a compatible triple $(J,\omega,g)$ for which each of the three structures is invariant.

The main result of this section is the following:

\begin{theorem}\label{thm:solvmanifolds}
    Let $M$ be a $4$-dimensional solvmanifold endowed with an invariant compatible triple $(J,\omega,g)$. Then the numbers $h^k_{d+d^c}$, $h^k_{d + d^\Lambda}$ are independent of the choice of invariant compatible triple.
\end{theorem}

The proof of theorem \ref{thm:solvmanifolds} follows from several lemmas that are valid in a slightly more general setting and that will be useful later for explicit  (see example \ref{ex:hopf}). Moreover, we will explicitly compute the numbers $h^k_{d+d^c}$, $h^k_{d + d^\Lambda}$ of compact complex surfaces diffeomorphic to solvmanifolds. The resulting numbers are summarized in tables \ref{tab:abc} and \ref{tab:def}.\\
We begin by taking care of $1$-forms.

\begin{lemma}\label{lemma:h1ddlambda}
    Let $M = \Gamma \backslash G$ be a compact quotient of a Lie group by a lattice endowed with an invariant almost symplectic structure $\omega$. Then $h^1_{d + d^\Lambda} = b_1$ and it is metric independent.
\end{lemma}
\begin{proof}
    Fix a compatible metric $g$. The inequality $h^1_{d + d^\Lambda} \le b_1$ holds for arbitrary almost Hermitian manifolds. For the opposite inequality, let $\alpha \in \H^1_d$ be a $d$-harmonic $1$-form. Then $d \alpha =0$ and $(d d^\Lambda)^* \alpha = d^c d^* \alpha =0$ for any choice of compatible metric. Since $\omega$ is invariant, $ *_s \alpha$ is an invariant $3$-form, thus it is $d$-closed and we have $d^\Lambda \alpha =0$.
\end{proof}

The following lemma deals with the space $\H^3_{d+d^\Lambda}$.

\begin{lemma}\label{lemma:h3ddlambda}
    Let $(M,J,\omega,g)$ be a Hermitian $2m$-manifold. Suppose that for every $\gamma \in \H^1_{ d d^c}$ we have
    \[
    (d^c)^* d^* d \gamma =0.
    \]
    Then $\H^{2m-1}_{d + d^\Lambda} = \H^{2m-1}_{d + d^c} \cap \ker \Delta_{d + d^\Lambda}$ and $h^{2m-1}_{d+d^\Lambda} \le h^{2m-1}_{d + d^c}$.
\end{lemma}
\begin{proof}
    By proposition \ref{prop:harmonic:iso}, $\H^{2m-1}_{d + d^\Lambda} \cong \H^1_{d^\Lambda d}$, thus we compute $d^\Lambda d$-harmonic $1$-forms. Let $\{ \gamma_1, \dots, \gamma_t \}$ be a basis of $\H^1_{d + d^c}$, let $\alpha \in \H^1_{d^\Lambda d}$ and let 
    \begin{equation}\label{eq:ddc:decomposition}
    \alpha = \sum_{j=1}^t A_j \gamma_j + (d d^c)^* \beta + df + d^c g
    \end{equation}
    be the harmonic decomposition of $\alpha$ with respect to $\Delta_{dd^c}$, with $A_j \in \C$, $\beta \in A^3$, $f,g \in C^\infty (M)$. Since $d^c \alpha =0$ and $\gamma_j \in \H^1_{dd^c}$, we have
    \[
    0 = d d^c \alpha = d d^c (d d^c)^* \beta,
    \]
    obtaining $(dd^c)^*\beta=0$. On the other side, we have that
    \[
    0 = d^\Lambda d \alpha = \sum_{j=1}^t A_j (d^c)^* d \gamma_j + (d^c)^* d d^c g.
    \]
    Taking the inner product with $dg$ and using the assumption $(d^c)^* d^* d \gamma_j =0$, we can write
    \[
    0 = \sum_{j=1}^t A_j \langle (d^c)^* d^* d \gamma_j, g \rangle + \langle d d^c g, d^c d g \rangle = - \norm{dd^c g}^2,
    \]
    which implies that $g$ is constant. Finally, from
    \[
    0 = d^* \alpha = d^* df
    \]
    we deduce that also $f$ is constant, $\alpha = \sum_j A_j \gamma_j \in \H^1_{d d^c}$ and $\H^1_{d^\Lambda d}$ is computed as $\H^1_{d d^c} \cap \ker \Delta_{d^\Lambda d}$.
\end{proof}

Denote by $\H^{p,q}_d$ the space of $d$-harmonic $(p,q)$-forms.

\begin{proposition}\label{prop:h3ddlambda}
    Let $M = \Gamma \backslash G$ be a $4$-dimensional compact quotient of a Lie group by a lattice and let $(J,\omega,g)$ be an invariant compatible triple, with $J$ integrable. Then
    \[
    \H^3_{d + d^\Lambda} = \H^{2,1}_d \cup \H^{1,2}_d.
    \]
\end{proposition}
\begin{proof}
    Since the compatible triple is invariant, any $1$-form $\gamma$ automatically satisfies $(d^c)^* d^* d \gamma =0$. By lemma \ref{lemma:h3ddlambda}, we have that
    \[
    \H^3_{d+d^\Lambda} \cong \H^1_{d^\Lambda d} = \H^1_{d d^c} \cap \ker \Delta_{d^\Lambda d} = A^1 \cap \ker \Delta_{d d^c} \cap \ker \Delta_{d^\Lambda d}.
    \]
    Let $\gamma \in A^1 \cap \ker \Delta_{d d^c} \cap \ker \Delta_{d^\Lambda d}$. Since $\Delta_{d d^c} \gamma =0$ and $\Delta_{d^\Lambda d} \gamma =0$, $\gamma$ must satisfy the following equations:
    \[
    d^c \gamma =0, \quad d^* \gamma =0, \quad (d^c)^* \gamma =0, \quad d^cd \gamma =0, \quad d^\Lambda d \gamma =0.
    \]
    We show that $d \gamma =0$. Let
    \begin{equation}\label{eq:dgamma}
    d \gamma = f \omega + \beta^{1,1} + \beta^{(2,0)(0,2)}
    \end{equation}
    be the Lefschetz and bidegree decomposition of $d \gamma$, with $f \in C^\infty (M)$ and $\beta^{1,1}$, $\beta^{(2,0)(0,2)}$ primitive forms. From the equations
    \begin{align*}
    0 &= d^2 \gamma = d ( f \omega + \beta^{1,1} + \beta^{(2,0)(0,2)} ), \\
    0 &= d^c d\gamma = - J d ( f \omega + \beta^{1,1} - \beta^{(2,0)(0,2)} ), \\
    0 &= d^\Lambda d\gamma = - *_s d ( f \omega - \beta^{1,1} - \beta^{(2,0)(0,2)} ),
    \end{align*}
    we get that $d ( f \omega) =0$, $d \beta^{1,1} =0$ and $d \beta^{(2,0)(0,2)} =0$. Since each of the terms in \eqref{eq:dgamma} is $d$-closed, we have that
    \[
    d^* d\gamma = - * d ( f \omega - \beta^{1,1} + \beta^{(2,0)(0,2)} ) =0,
    \]
    which implies $d \gamma =0$ and that $\gamma$ is both $d$-harmonic and $d^c$-harmonic. As a consequence, the $(1,0)$-degree part and $(0,1)$-degree part of $\gamma$ are both $d$-harmonic. This shows that $\H^1_{d^\Lambda d} = \H^{1,0}_d \cup \H^{0,1}_d$ and, after applying the isomorphism $*$, it concludes the proof of the proposition.
\end{proof}

We are ready for the proof of theorem \ref{thm:solvmanifolds}.

\begin{proof}[Proof of theorem \ref{thm:solvmanifolds}.]
Let $M$ be a $4$-dimensional solvmanifold admitting a complex structure. Then $M$ is one of the following \cite{Has05}:
\begin{itemize}
    \item [(A)] complex torus;

    \item [(B)] hyperelliptic surface;

    \item [(C)] Inoue surface of type $\mathcal{S}_M$;

    \item [(D)] primary Kodaira surface;

    \item [(E)] secondary Kodaira surface;

    \item [(F)] Inoue surface of type $\mathcal{S}^\pm$.
\end{itemize}
By theorem \ref{thm:BC:diamond}, the numbers $h^k_{d+d^c}$ depend only on $b_1$, $b^+$, $b^-$ and not on the choice of compatible triple (not necessarily invariant).\\
All invariant structures on the torus and the hyperelliptic surface (cases (A) and (B)) are K\"ahler structures, therefore $h^k_{d+d^\Lambda} = b_k$ and they are independent of the compatible triple.\\
Cases (C), (E) and (F) can be treated simultaneously since they have the same Betti numbers $b_1=1$, $b_2=0$. By theorem \ref{thm:BC:diamond} we have that $h^1_{d+d^c} = b_1-1 = 0$, $h^2_{d+d^c}= b_2 +1 =1$, and $h^3_{d+d^c}= b_1+1 = 2$. By lemma \ref{lemma:h1ddlambda}, $h^1_{d + d^\Lambda} = b_1 = 1$. By theorem \ref{thm:hermitian:betti} and the fact that there are no $d$-harmonic $2$-forms since $b_2=0$, we have that $h^2_{d + d^\Lambda} = b_2 =0$. By proposition \ref{prop:h3ddlambda}, we have that $\H^3_{d + d^\Lambda} = \H^{2,1}_d \cup \H^{1,2}_d$. In particular, since $\H^{2,1}_d \cong \H^{1,2}_d$ via complex conjugation, $h^3_{d+d^\Lambda}$ must be even. Moreover, $\H^{2,1}_d \cap \H^{1,2}_d = \{ 0 \}$ and both spaces inject into $\H^3_d$, thus $h^3_{d+d^\Lambda} \le b_1 = 1$ and it must be an even number, showing that $h^3_{d+d^\Lambda}=0$.\\
Case (D) has to be treated separately, since the Betti numbers in this case are $b_1 =3$, $b_2 =4$. By theorem \ref{thm:BC:diamond}, we immediately have $h^1_{d+d^c} = 2$, $h^2_{d+d^c} = 5$, and $h^1_{d+d^c} = 4$. By lemma \ref{lemma:h1ddlambda}, we have that $h^1_{d+d^\Lambda} = b_1 = 3$. By proposition \ref{prop:h3ddlambda}, with the same reasoning as in cases (C), (E) and (F), we know that $h^3_{d+d^\Lambda}$ is even and $h^3_{d+d^\Lambda} \le b_1 = 3$, so that either $h^3_{d+d^\Lambda} =0$ or $h^3_{d+d^\Lambda}=2$. We show that $h^3_{d+d^\Lambda}=2$ by showing that there is at least one $(d+d^\Lambda)$-harmonic $3$-form. Fix an invariant compatible triple $(J,\omega,g)$. By \cite[Theorem 1.3]{Sal01}, since primary Kodaira surfaces are nilmanifolds $\Gamma \backslash G$ and the complex structure is invariant, up to a linear transformation of the Lie algebra of $G$, we can find a basis of invariant $(1,0)$-form $\{\phi^1, \phi^2 \}$ such that $d \phi^1 =0$, $d \phi^2 = \phi^{1 \bar 1}$. Since $g$ is an invariant metric, $*\phi^1$ is an invariant $(2,1)$-form, thus $d( * \phi^1)=0$. Moreover, $d^*(* \phi^1) = - * d \phi^1 =0$. Hence, $* \phi^1 \in \H^{2,1}_d$ and $h^3_{d+d^\Lambda} =2$.
To conclude the proof, we show that $\H^2_d = \H^2_{d^c} \subseteq \H^2_{d+d^\Lambda}$, which implies that $h^2_{d+d^\Lambda} =b_2 = 4$ by theorem \ref{thm:hermitian:betti}. Consider the decomposition $\H^2_d = \H^+_g \oplus \H^-_g$. Forms in $\H^-_g$ are anti-self-dual, thus they have bidegree $(1,1)$ and they are primitive. Moreover, they are also $d$-closed. If $\alpha \in \H^-_g$, then $d \alpha =0$, $d^\Lambda \alpha = - *_s d *_s \alpha = *_s d \alpha =0$ and $(d d^\Lambda)^* \alpha = d^c d^* \alpha = d^c * d \alpha =0$. As a consequence, we have $\H^-_g \subseteq \H^2_{d+d^\Lambda}$. For the inclusion $\H^+_g \subseteq \H^2_{d+d^\Lambda}$, we observe that on primary Kodaira surfaces $b^+=2$ and, as long as the compatible triple is invariant, we have $\H^+_g = \C \langle \phi^{12}, \phi^{\bar 1 \bar 2} \rangle$, where $\{ \phi^1, \phi^2 \}$ is the preferred basis of \cite{Sal01} we considered above in the proof, up to normalization. Indeed, we have that $d \phi^{12} =0$ and $* \phi^{12}$ is an invariant $(2,0)$-form since the metric is invariant. Finally, also $*_s \phi^{12}$ is an invariant $(2,0)$-form since $\omega$ is invariant, and one has that $d \phi^{12} =0$, $d^\Lambda \phi^{12} = - *_s d *_s \phi^{12} =0$ and $(d d^\Lambda)^* \phi^{12} = d^c d^* \phi^{12} =0$, proving that $\H^+_g \subseteq \H^2_{d+d^\Lambda}$ and $h^2_{d+d^\Lambda} =4$.
\end{proof}
    
The results proved in this section and in section \ref{sec:surfaces} can be used to compute $h^k_{d + d^c}$, $h^k_{d + d^\Lambda}$ in the case of complex surfaces not necessarily diffeomorphic to solvmanifolds, as we illustrate in the example below.

\begin{example}[Hopf surface]\label{ex:hopf}
    Let $M= \S^1 \times \S^3$ be the Hopf surface. There exists a parallelism for $T^*M$ given by $\{ e^1, e^2, e^3, e^4 \}$, with differentials
    \[
    d e^1 =0, \quad d e^2 = - e^{34}, \quad d e^3 = e^{24}, \quad d e^4 = - e^{23}.
    \]
    Note that $M$ is the quotient of a non-solvable Lie group by a lattice.

    \begin{proposition}\label{prop:hopf}
        Let $(M,J,\omega,g)$ be the Hopf surface endowed with an invariant compatible triple with $J$ integrable. The numbers $h^k_{d + d^c}$, $h^k_{d + d^\Lambda}$ are 
        \begin{table}[!h]
            \centering
            \begin{tabular}{c|c c c c c}
                 $k$ & $0$ & $1$ & $2$ & $3$ & $4$  \\
                 \hline
                 $h^k_{d+d^c}$ & $1$ & $0$ & $1$ & $2$ & $1$ \\
                 $h^k_{d+d^\Lambda}$ & $1$ & $1$ & $0$ & $0$ & $1$
            \end{tabular}
        \end{table}
    \end{proposition}
    \begin{proof}
        By theorem \ref{thm:BC:diamond}, we have that $h^1_{d +d^c} = 0$, $h^2_{d + d^c} = 1$ and $h^3_{d + d^c} = 2$.\\
        To compute $h^k_{d +d^\Lambda}$, we resort to several different arguments. By lemma \ref{lemma:h1ddlambda}, $h^1_{d + d^\Lambda} = b_1 =1$. Since there is no $d$-harmonic $2$-form, theorem \ref{thm:hermitian:betti} tells us that $h^2_{d+d^\Lambda} = 0$. Finally, by proposition \ref{prop:h3ddlambda} and the same reasoning used in the proof of theorem \ref{thm:solvmanifolds}, we conclude that $h^3_{d +d^\Lambda} = 0$.
    \end{proof}
\end{example}

\begin{remark}
    The fact that almost complex and almost symplectic invariants are determined by the topology of the underlying manifold is not surprising, especially on solvmanifolds (see, e.g., \cite[Theorem 5.15]{ST23a}).
\end{remark}

    \begin{table}
    \caption{\label{tab:abc}The numbers $h^k_{d+d^c}$, $h^k_{d+d^\Lambda}$ of the complex torus, the hyperelliptic surface and the Inoue surface $\mathcal{S}_M$.}
    \begin{adjustbox}{center}
    \begin{tabular}{ c|c c|c c|c c}
      & \multicolumn{2}{c|}{(A) Complex torus} & \multicolumn{2}{c|}{\text{(B) Hyperelliptic surface}} & \multicolumn{2}{c}{\text{(C) Inoue surface $\mathcal{S}_M$}} \\ \cline{2-7}
     & & & & & & \\[-2.3ex]
    $k$ & $h^k_{d+d^c}$ & $h^k_{d+d^\Lambda}$ & $h^k_{d+d^c}$ & $h^k_{d+d^\Lambda}$ & $h^k_{d+d^c}$ & $h^k_{d+d^\Lambda}$ \\
    \hline
    0 & 1 & 1 & 1 & 1 & 1 & 1 \\
    1 & 4 & 4 & 2 & 2 & 0 & 1 \\
    2 & 6 & 6 & 2 & 2 & 1 & 0 \\
    3 & 4 & 4 & 2 & 2 & 2 & 0 \\
    4 & 1 & 1 & 1 & 1 & 1 & 1 \\
    \end{tabular}
    \end{adjustbox}
    \end{table}

    \begin{table}
    \caption{\label{tab:def}The numbers $h^k_{d+d^c}$, $h^k_{d+d^\Lambda}$ of the primary and secondary Kodaira surface and the Inoue surface $\mathcal{S}^\pm$.}
    \begin{adjustbox}{center}
    \begin{tabular}{ c|c c|c c|c c}
      & \multicolumn{2}{c|}{(D) Primary Kodaira surface} & \multicolumn{2}{c|}{\text{(E) Secondary Kodaira surface}} & \multicolumn{2}{c}{\text{(F) Inoue surface $\mathcal{S}^\pm$}} \\ \cline{2-7}
     & & & & & & \\[-2.3ex]
    $k$ & $h^k_{d+d^c}$ & $h^k_{d+d^\Lambda}$ & $h^k_{d+d^c}$ & $h^k_{d+d^\Lambda}$ & $h^k_{d+d^c}$ & $h^k_{d+d^\Lambda}$ \\
    \hline
    0 & 1 & 1 & 1 & 1 & 1 & 1 \\
    1 & 2 & 3 & 0 & 1 & 0 & 1 \\
    2 & 5 & 4 & 1 & 0 & 1 & 0 \\
    3 & 4 & 2 & 2 & 0 & 2 & 0 \\
    4 & 1 & 1 & 1 & 1 & 1 & 1 \\
    \end{tabular}
    \end{adjustbox}
    \end{table}

{\small
\printbibliography
}

\end{document}